\documentclass[12pt, a4paper]{article}
\usepackage{amsmath}
\usepackage{comment}
\usepackage{amsfonts}
\usepackage{amssymb}
\usepackage{epsfig}
\usepackage{graphicx}
\usepackage{color}
\usepackage{amsthm}
\usepackage{enumerate}
\usepackage [latin1]{inputenc}
\usepackage[numbers, sort&compress]{natbib}
\usepackage{url}
\usepackage{subcaption}
\newtheorem{thm}{Theorem}[section]
\newtheorem{lem}[thm]{Lemma}

\theoremstyle{definition}

\graphicspath{{figures/}}

\usepackage{url}

\usepackage{authblk}

\long\def\delete#1{}

\usepackage{xcolor}
\usepackage[normalem]{ulem}

\begin{document}
\openup 0.5\jot
\underline{}\begin{center}{\Large \bf Maximum Forcing Number of $TUC_4C_8$ Nanotubes and Nanorings}\\
\vspace{4mm}
{Shouang Sun, Qianqian Liu\footnote{The corresponding author. \textsf{E-mail: liuqq2023@imut.edu.cn}}}
\vspace{2mm}

{ College of Science, Inner Mongolia University of Technology, Hohhot,\\ Inner Mongolia 010010, P. R. China}\\
\end{center}

\setlength{\baselineskip}{20pt}
\noindent {\bf Abstract}:
The maximum forcing number is an important topological index reflecting the chemical stability of molecules. Carbon nanotubes, distinguished by their unique polycyclic configurations, are valuable in chemical synthesis, biosensing, and the design of novel functional materials. Diudea et al. (2004) constructed two classes of 3-regular graphs on carbon nanotubes, $TUC_4C_8(R)$ and $TUC_4C_8(S)$ nanotubes, denoted by $T_{p,q}$ and $N_{p,q}$, where \(p\) and \(q\) count the rhombi in the horizontal rows and vertical columns of \(T_{p,q}\), and the regular octagons in those of \(N_{p,q}\), respectively. When the vertices on the two horizontal boundaries are joined pairwise, these nanotubes close up into the corresponding nanorings $\mathrm{Tori}_{p,q}$ and $H_{p,q}$. Jiang and Zhang (2016) determined the maximum forcing number of $\mathrm{Tori}_{p,q}$ to be $pq$.  In this paper we determine all remaining  cases: $T_{p,q}$ and $H_{p,q}$ have maximum forcing numbers $pq$ and $2pq$, respectively, while $N_{p,q}$ has maximum forcing number $2q$ when $p=1$ and $p(2q-1)$ when $p>1$.

\vspace{2mm} \noindent{\textbf{Keywords:}}  Perfect matching, maximum forcing number, nanotube, nanoring
\vspace{2mm}


\section{Introduction}
In chemical graph theory, molecular structures can be naturally abstracted into graph-theoretic models: atoms correspond to the vertices of a graph, and the chemical bonds between atoms correspond to the edges joining these vertices. In this way, a Kekul\'{e} structure of a molecule is represented by a perfect
matching of its graph. Hence the analysis of resonance structures and molecular stability reduces to the quantitative study of perfect matchings in
graph theory.

In 1985, Klein and Randi\'{c} \cite{klein87} observed that fixing some double bonds of a Kekul\'{e} structure determines all the remaining ones, and introduced the internal degree of freedom, the minimum number of double bonds needed to determine a Kekul\'{e} structure. This notion was later named the forcing number by Harary, Klein and \v{Z}ivkovi\'{c} \cite{harary91}.

Let $G$ be a simple graph with the vertex set $V(G)$ and edge set $E(G)$. Let $M$ be a perfect matching of $G$. A subset $S\subseteq M$ is called a \emph{forcing set} of $M$ if $S$ is not contained in any other perfect matching of $G$. The minimum number of edges in a forcing set of $M$ is called the \emph{forcing number} of $M$, denoted by $f(G,M)$. The \emph{maximum forcing number} of $G$, denoted by $F(G)$, is the maximum value of $f(G,M)$ over all perfect matchings $M$ of $G$. From the perspective of structural chemistry, the maximum forcing number of a
molecular graph not only quantifies the largest possible constraint on the set
of perfect matchings but also reflects the topological rigidity of the
molecular skeleton. However, the computational complexity of determining the
maximum forcing number of a graph is still open~\cite{afshani04}. We refer the reader to the survey of Zhang et al. \cite{zhang25} on problems related to forcing.

Since the microscopic properties of real nanotubes are difficult to measure experimentally, chemical graph theory has become a key bridge between molecular structure and macroscopic properties. Motivated by this, Diudea and his co-workers~\cite{diudea05,stefu04}
constructed the $C_4C_8$ grid, a cubic planar lattice obtained by arranging
squares and octagons alternately. Wrapping this grid into a cylinder in two
different ways yields two distinct nanotubes, $TUC_4C_8(R)$ and
$TUC_4C_8(S)$. Identifying the corresponding vertices of the two horizontal
boundaries in pairs then closes each of them into a nanoring.
Various topological indices, such as Wiener index \cite{stefu04}, Hosoya polynomial \cite{xuzhang09}, and PI index \cite{xudeng07,iranmanesh07}, have been computed for these nanotubes and their corresponding nanorings.
Jiang and Zhang \cite{jiangzhang16} determined the maximum forcing number of $TUC_4C_8(R)$ nanoring, but no results are available for the other nanotubes and nanorings.

In this paper, we derive explicit expressions for the maximum forcing numbers of $TUC_4C_8(R)$ and $TUC_4C_8(S)$ nanotubes and  that of $TUC_4C_8(S)$ nanorings, with a view to providing theoretical support for assessing the structural stability of nanomaterials at the microscopic topological level. In Section 2, we present some notations and two key lemmas. In Section 3, we determine the maximum forcing numbers of $TUC_4C_8(R)$ nanotubes; in Section 4, we determine those of $TUC_4C_8(S)$ nanotubes and nanorings.
\section{\normalsize Preliminaries}%
In this section, we introduce some notations and state two lemmas that will be used repeatedly in Sections 3 and 4.

For a positive integer $p$ and each $0\leq i\leq p-1$, define
\[
\mathbb{Z}_{p,i}=\{i,i+1,\ldots,p-1\},
\]
so that $\mathbb{Z}_{p,0}=\mathbb{Z}_p=\{0,1,2,\ldots,p-1\}$. In addition, set $\mathbb{Z}_{p,-1}=\{0,1,2,\ldots,p-2\}$.

For a subset $T\subseteq V(G)$, the subgraph of $G$ \emph{induced by} $T$, denoted by $G[T]$, has vertex set $T$ and edge set $\{xy\in E(G) \mid x,y\in T\}$. For a subset $S\subseteq E(G)$, let $V(S)$ denote the set of vertices incident with the edges in $S$.  We write $G-V(S)$ for the subgraph of $G$ induced by
$V(G)\setminus V(S)$, that is, the graph obtained from $G$ by deleting all vertices in $V(S)$ together with their incident edges.
A \emph{$2$-path} in $G$ is an induced path $uvw$ of length two, whose vertices
$u$ and $w$ are called \emph{endvertices}. When we say that two cycles in $G$ are \emph{disjoint}, we mean that their vertex sets are disjoint.

Let $G$ be a graph with a perfect matching $M$, and let $C$ be a cycle of even length in $G$. If the edges of $C$ appear alternately in $M$ and in $E(G)\setminus M$, then we call $C$ is $M$-\emph{alternating}. Adams et al. \cite{adams04} and Kleinerman \cite{kleinerman06} proved the following results, independently.
\begin{lem}\cite{adams04,kleinerman06}\label{lem:alternating}
Let $M$ be a perfect matching of $G$ and let $S\subseteq M$. The following statements are equivalent:
\begin{enumerate}
\renewcommand{\labelenumi}{(\roman{enumi})}
\item $S$ is a forcing set of $M$;
\item every $M$-alternating cycle of $G$ contains at least  an edge of $S$;
\item $G-V(S)$ contains no $M$-alternating cycles.
\end{enumerate}
\end{lem}

To obtain upper bounds for the maximum forcing number, it is useful to delete a set
of vertices that can lie on only few $M$-alternating cycles. Given an independent
set $T$ of $G$, its vertices are called \emph{marked vertices}. Define
\[
M_T=\{e\in M \mid e \text{ has an endvertex in } T\}.
\]
The following lemma shows
that if the union of all 2-paths with both endvertices in $T$ carries no
$M$-alternating cycle, then $M\setminus M_T$ is a forcing set of $M$.

\begin{lem}\label{lem:labeled}
Let $G$ be a graph with a perfect matching $M$, and $T$ be an independent set
of $G$. Let $G'$ be the union of all 2-paths of $G$ whose endvertices lie in $T$. If $G'$ contains no
$M$-alternating cycle, then $f(G,M)\le |M|-|T|$.
\end{lem}
\begin{proof}We prove that $G[V(M_T)]$ contains no $M$-alternating cycles.
Suppose to the contrary that $G[V(M_T)]$ contains an $M$-alternating cycle $C$. Then $C$ is also an $M_T$-alternating cycle. Since $M$ is a perfect matching and $T$ is an independent set, each edge of $M_T$ contains exactly one vertex of $T$, and each vertex of $T$ lies in exactly one edge of $M_T$. Thus the map sending every vertex of $T$ to the unique edge of $M$ incident with it is a bijection from $T$ onto $M_T$, and $|M_T|=|T|$. Furthermore, half vertices of $C$ are marked, and marked and unmarked vertices appear alternately along $C$. Thus, $C$ is the union of 2-paths of $G$ whose endvertices lie in $T$. Hence $C$ is an $M$-alternating cycle contained in $G'$, which is a contradiction.

By Lemma \ref{lem:alternating}, $M\backslash M_T$ is a forcing set of $M$. Thus,
\[
f(G,M)\leq |M\setminus M_T|=|M|-|T|.
 \qedhere
\]
\end{proof}

\section{$TUC_4C_8(R)$ Nanotubes}\label{sec:tuc4c8r}
\begin{figure}[htbp]
  \centering
  \includegraphics[height=5cm,width=9.5cm]{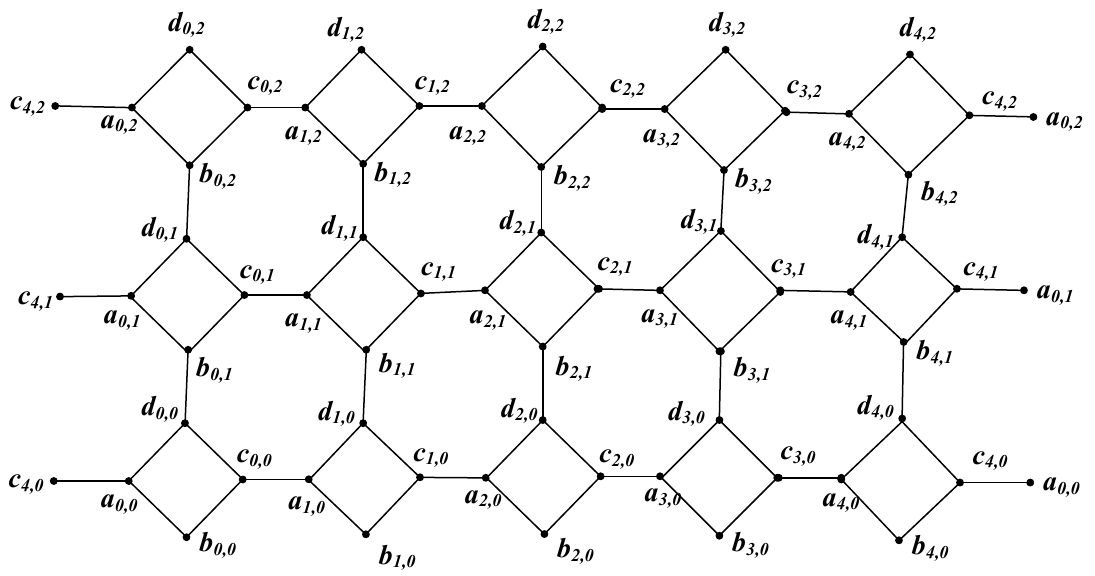}
  \caption{The rectangular representation of nanotube $T_{5,3}$ with the labeling of its vertices.}
  \label{fig:r-rectangle}
\end{figure}

The $TUC_4C_8(R)$ \emph{nanotube} is a portion of a cylindrical tiling in which rhombi and regular octagons alternate. This nanotube, obtained by gluing the two vertical boundaries of a rectangle carrying the $C_4C_8$ lattice (see Fig.~\ref{fig:r-rectangle}), is denoted by $T_{p,q}$, where $p$ and $q$ are the numbers of rhombi in each horizontal row and each vertical column of the rectangle, respectively.
Two rhombi are said to be \emph{in the same row} if the line through their centers is parallel to the horizontal sides, and \emph{in the same column} if it is parallel to the vertical sides.

We introduce the following coordinate system. The rhombi of $T_{p,q}$ form a $p\times q$ array, whose rows are numbered $0,1,\ldots,q-1$ from bottom to top and whose columns are numbered $0,1,\ldots,p-1$ from left to right. Hence the rhombus in column $i$ and row $j$ is denoted by $r_{i,j}$, where $i\in\mathbb{Z}_p$ and $j\in\mathbb{Z}_q$, and its four vertices are labeled counterclockwise by $a_{i,j}$, $b_{i,j}$, $c_{i,j}$, $d_{i,j}$, starting from the leftmost vertex (see Fig. \ref{fig:r-rectangle}).
When the vertices on the two horizontal boundaries are joined pairwise, the $TUC_4C_8(R)$  nanotubes $T_{p,q}$ close up into the corresponding \emph{nanorings} $\mathrm{Tori}_{p,q}$.

\begin{thm}\label{thm:main} Let $p\ge 1$ and $q\ge 1$. Then the maximum forcing number  $F(T_{p,q})=pq.$
\end{thm}

\begin{proof}
We first establish the lower bound $F(T_{p,q})\ge pq$.  Consider the perfect matching
\[
M_0=\{a_{i,j}d_{i,j},\ b_{i,j}c_{i,j}\mid i\in\mathbb{Z}_p,\ j\in\mathbb{Z}_q\}
\]
(see red edges in Fig.~\ref{rfig2}).
\begin{figure}[htbp]
  \centering
  \begin{minipage}{0.48\textwidth}
    \centering
    \includegraphics[width=\linewidth]{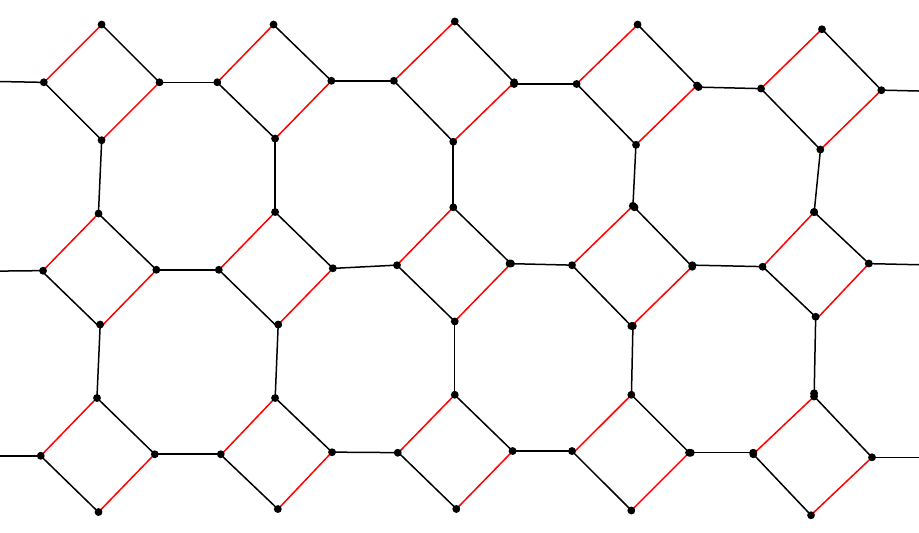}
    \caption{The perfect matching $M_0$ of $T_{5,3}$.}
    \label{rfig2}
  \end{minipage}\hfill
  \begin{minipage}{0.48\textwidth}
    \centering
    \includegraphics[width=\linewidth]{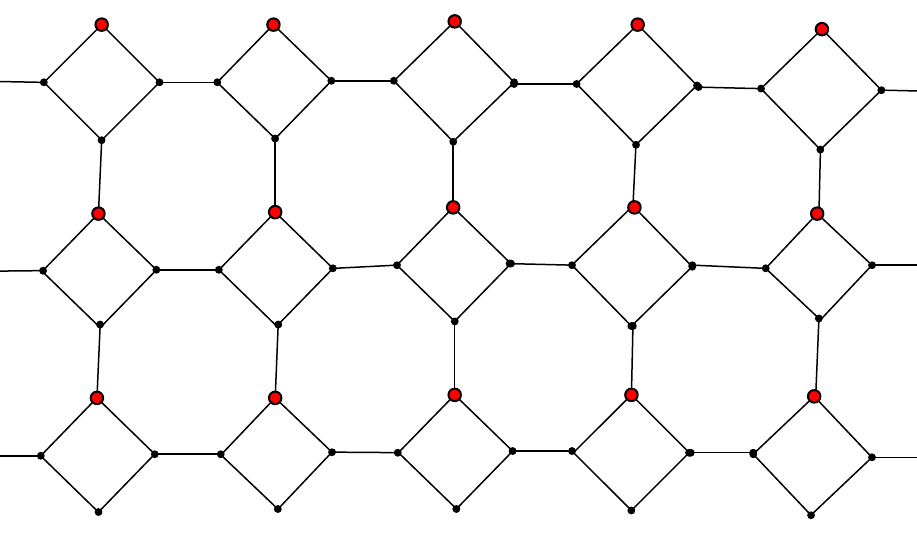}
    \caption{The marked vertices $T$ of $T_{5,3}$.}
    \label{rfig3}
  \end{minipage}
\end{figure}
We show that $f(T_{p,q},M_0)=pq$. Let
\[
S_0=\{a_{i,j}d_{i,j}\mid i\in\mathbb{Z}_p,\ j\in\mathbb{Z}_q\}\subseteq M_0.
\]
Since $T_{p,q}$ is cubic, the neighbors of $b_{i,j}$ lying outside $V(S_0)$
are precisely $c_{i,j}$, and the neighbors of $c_{i,j}$ lying outside
$V(S_0)$ are precisely $b_{i,j}$, for all $i\in\mathbb{Z}_p$ and
$j\in\mathbb{Z}_q$; indeed, the remaining neighbors $a_{i,j}$, $d_{i,j-1}$ of
$b_{i,j}$ and $d_{i,j}$, $a_{i+1,j}$ of $c_{i,j}$ (whenever these vertices
exist) all lie in $V(S_0)$. Hence $T_{p,q}-V(S_0)$ consists of the $pq$ edges
$b_{i,j}c_{i,j}$ with $i\in\mathbb{Z}_p$ and $j\in\mathbb{Z}_q$, and these
edges form its unique perfect matching. 
Thus, $S_0$ is a forcing set of $M_0$, and
$f(T_{p,q},M_0)\le |S_0|=pq$. On the other hand, every rhombus of $T_{p,q}$ is an
$M_0$-alternating cycle, and these $pq$ rhombi are pairwise disjoint. Since a
forcing set of $M_0$ must meet every $M_0$-alternating cycle, $f(T_{p,q},M_0)\ge pq$ by Lemma \ref{lem:alternating}. Therefore
\[
F(T_{p,q})\ge f(T_{p,q},M_0)=pq.
\]

For the reverse inequality, let $M$ be an arbitrary perfect matching of $T_{p,q}$. We take the marked vertices
$T=\{d_{i,j}\mid i\in\mathbb{Z}_p,\ j\in\mathbb{Z}_q\}$
(see red vertices in Fig.~\ref{rfig3}). Since the distance between any two distinct vertices of $T$ is at least $3$, no
2-path of $T_{p,q}$ has both endvertices in $T$. Therefore the graph $G'$ formed by
the union of all such 2-paths has no edges, and it contains no $M$-alternating cycle. By Lemma \ref{lem:labeled},
\[
f(T_{p,q},M)\le |M|-|T|=2pq-pq=pq.
\]
Since $f(T_{p,q},M)\le pq$ for every perfect matching $M$ of $T_{p,q}$, and
$F(T_{p,q})$ is the maximum of these values, we have $F(T_{p,q})\le pq$. Together
with the reverse inequality $F(T_{p,q})\ge pq$ proved above, we obtain
\[
F(T_{p,q})=pq. \qedhere
\]
\end{proof}

\section{$TUC_4C_8(S)$ Nanotubes and nanorings}\label{sec:s}
\begin{figure}[htbp]
\centering
\includegraphics[width=0.70\textwidth]{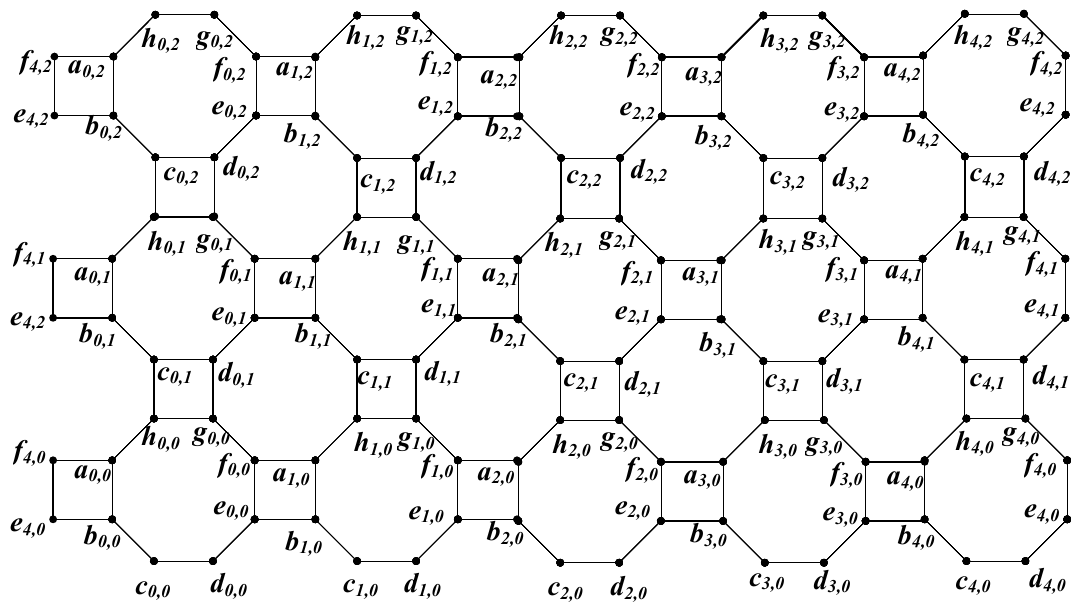}
\caption{The rectangular representation of nanotube $N_{5,3}$ with the labeling of its vertices.}
\label{sfig1}
\end{figure}

The $TUC_4C_8(S)$ \emph{nanotube} is a portion of a cylindrical tiling in which squares and regular octagons alternate. This nanotube, obtained by gluing the two vertical boundaries of a rectangle carrying the $C_4C_8$ lattice (see Fig.~\ref{sfig1}), is denoted by $N_{p,q}$, where $p$ and $q$ are the numbers of regular octagons in each horizontal row and each vertical column of the rectangle, respectively.
Furthermore, if we join the corresponding vertices on the two horizontal boundaries
of the rectangle in $N_{p,q}$ pairwise, we obtain the $TUC_4C_8(S)$ \emph{nanorings} $H_{p,q}$.  Two regular
octagons of $N_{p,q}$ and $H_{p,q}$ are said to be \emph{in the same row} if the line through their centers is
parallel to the horizontal sides, and \emph{in the same column} if it is parallel to the
vertical sides.

We introduce the following coordinate system. The regular octagons of $N_{p,q}$ and
$H_{p,q}$ form a $p\times q$ array, whose rows are numbered $0,1,\ldots,q-1$ from
bottom to top and whose columns are numbered $0,1,\ldots,p-1$ from left to right.
Hence the regular octagon in column $i$ and row $j$ is denoted by $O_{i,j}$, where
$i\in\mathbb{Z}_p$ and $j\in\mathbb{Z}_q$, and its eight vertices are labeled
counterclockwise by \(a_{i,j}\), \(b_{i,j}\), \(c_{i,j}\), \(d_{i,j}\),
\(e_{i,j}\), \(f_{i,j}\), \(g_{i,j}\), \(h_{i,j}\), starting from the upper-left vertex (see Fig.~\ref{sfig1}).

\begin{thm}\label{thm:s-tube}
Let $p\ge 1$ and $q\ge 1$. Then the maximum forcing number
\[
F(N_{p,q})=
\begin{cases}
2q, & p=1,\\[2pt]
p(2q-1), & p>1.
\end{cases}
\]
\end{thm}

\begin{proof}
Based on the value of $p$, we consider the cases $p=1$ and $p\ge2$.

For \(p=1\), we prove that $F(N_{1,q})=2q$ for all $q\ge 1$.
Consider the perfect matching
\[
M_0=\{g_{0,j}h_{0,j},\ a_{0,j}f_{0,j},\ b_{0,j}e_{0,j},\ c_{0,j}d_{0,j}\mid
j\in\mathbb{Z}_q\}
\]
of $N_{1,q}$, indicated by the red edges in Fig.~\ref{sfig3}. Let
\[
S_0=\{g_{0,j}h_{0,j},\ b_{0,j}e_{0,j}\mid j\in\mathbb{Z}_q\}
\]
be the subset of $M_0$. 
Since $N_{p,q}$ is cubic, the neighbors of $a_{0,j}$, $f_{0,j}$, $c_{0,j}$ and
$d_{0,j}$ lying outside $V(S_0)$ are precisely $f_{0,j}$, $a_{0,j}$, $d_{0,j}$
and $c_{0,j}$, respectively; indeed, the remaining neighbors $b_{0,j}$ and
$h_{0,j}$ of $a_{0,j}$, $e_{0,j}$ and $g_{0,j}$ of $f_{0,j}$, $b_{0,j}$ and
$h_{0,j-1}$ of $c_{0,j}$, and $e_{0,j}$ and $g_{0,j-1}$ of $d_{0,j}$ all
lie in $V(S_0)$ (whenever these vertices exist). Hence $N_{p,q}-V(S_0)$
consists of the $2q$ edges
$\{a_{0,j}f_{0,j},\ c_{0,j}d_{0,j}\mid j\in\mathbb{Z}_q\}$, and these edges
form its unique perfect matching. 
Thus $S_0$ is a forcing
set of $M_0$, and $f(N_{1,q},M_0)\le |S_0|=2q$. For the reverse inequality, \(N_{1,q}\) contains \(q\) pairs of \(M_0\)-alternating cycles
\[
A_j=a_{0,j}h_{0,j}g_{0,j}f_{0,j}a_{0,j},
\qquad
B_j=b_{0,j}c_{0,j}d_{0,j}e_{0,j}b_{0,j},
\qquad j\in\mathbb Z_q,
\]
indicated in green and blue, respectively, in Fig.~\ref{sfig8}. Since these \(2q\) cycles are pairwise disjoint, Lemma~\ref{lem:alternating} implies that every forcing set of \(M_0\) contains at least one edge of \(M_0\) from each \(M_0\)-alternating cycle. Since these \(2q\) edges are distinct, we have $f(N_{1,q},M_0)\ge 2q.$
Together with the upper bound \(f(N_{1,q},M_0)\le 2q\), this gives $f(N_{1,q},M_0)=2q.$
Therefore,
\[
F(N_{1,q})\ge f(N_{1,q},M_0)=2q.
\]
\begin{figure}[htbp]
  \centering
  \begin{minipage}{0.32\textwidth}
    \centering
    \includegraphics[height=7cm]{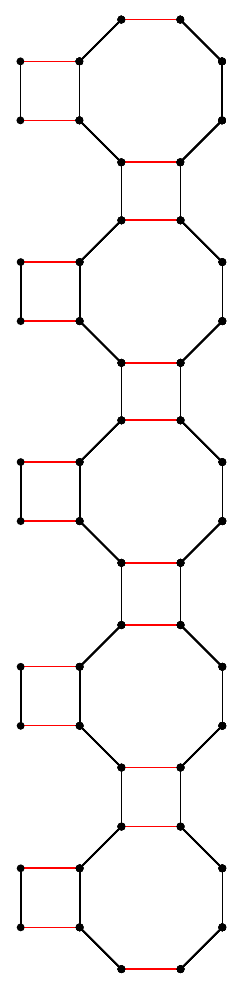}
    \caption{Perfect matching $M_0$.}
    \label{sfig3}
  \end{minipage}\hfill
  \begin{minipage}{0.32\textwidth}
    \centering
    \includegraphics[height=7cm]{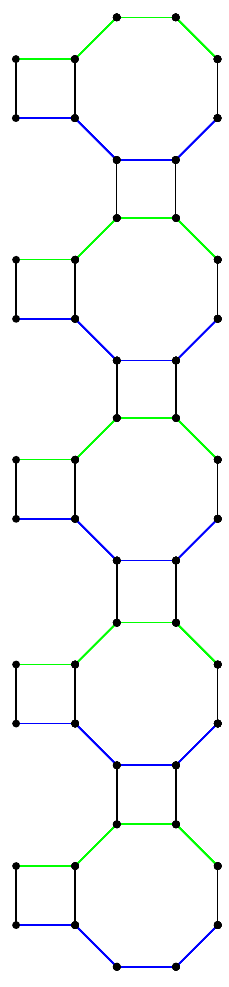}
    \caption{$M_0$-alternating cycles.}
    \label{sfig8}
  \end{minipage}\hfill
  \begin{minipage}{0.32\textwidth}
    \centering
    \includegraphics[height=7cm]{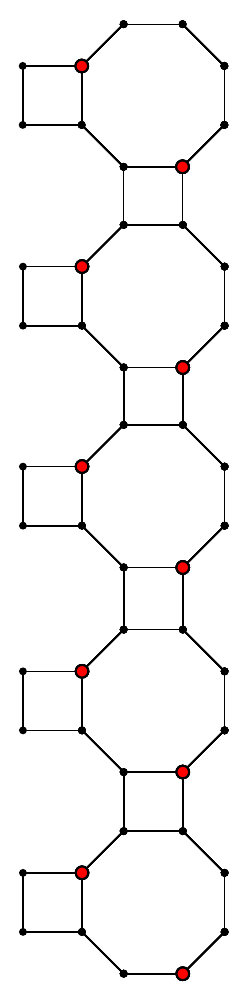}
    \caption{Marked vertices $T$.}
    \label{sfig4}
  \end{minipage}
\end{figure}

Conversely, we show that $F(N_{1,q})\le 2q$. Let $M$ be an arbitrary
perfect matching of $N_{1,q}$. We choose an independent set \(T=\{a_{0,j},d_{0,j}\mid j\in\mathbb{Z}_q\}\) of cardinality \(2q\), indicated by the red vertices in Fig.~\ref{sfig4}. Let $G'$ be the graph formed by the union of all 2-paths of $N_{1,q}$
whose endvertices lie in $T$. Since any two distinct vertices of \(T\) are at distance at least \(3\), no \(2\)-path joins two vertices of \(T\). Therefore \(G'\) has no edges and hence contains no \(M\)-alternating cycle.  By Lemma~\ref{lem:labeled},
\[
f(N_{1,q},M)\le |M|-|T|=4q-2q=2q.
\]
By the definition of the maximum forcing number, we have $F(N_{1,q})\le 2q$.
Combining this with the inequality \(F(N_{1,q})\ge 2q\) established above, we obtain \(F(N_{1,q})=2q\).

For \(p\ge 2\), we prove that \(F(N_{p,q})=p(2q-1)\) for all \(q\ge 1\).  We first establish the lower bound
$F(N_{p,q})\ge p(2q-1)$. Consider the perfect matching
\[
M_0=\{a_{i,j}b_{i,j},\ c_{i,j}d_{i,j},\ e_{i,j}f_{i,j},\ g_{i,j}h_{i,j}\mid
i\in\mathbb{Z}_p,\ j\in\mathbb{Z}_q\}
\]
of $N_{p,q}$ (see red edges in  Fig. \ref{sfig5}).
\begin{figure}[htbp]
  \centering
  \begin{minipage}{0.48\textwidth}
    \centering
    \includegraphics[width=\linewidth]{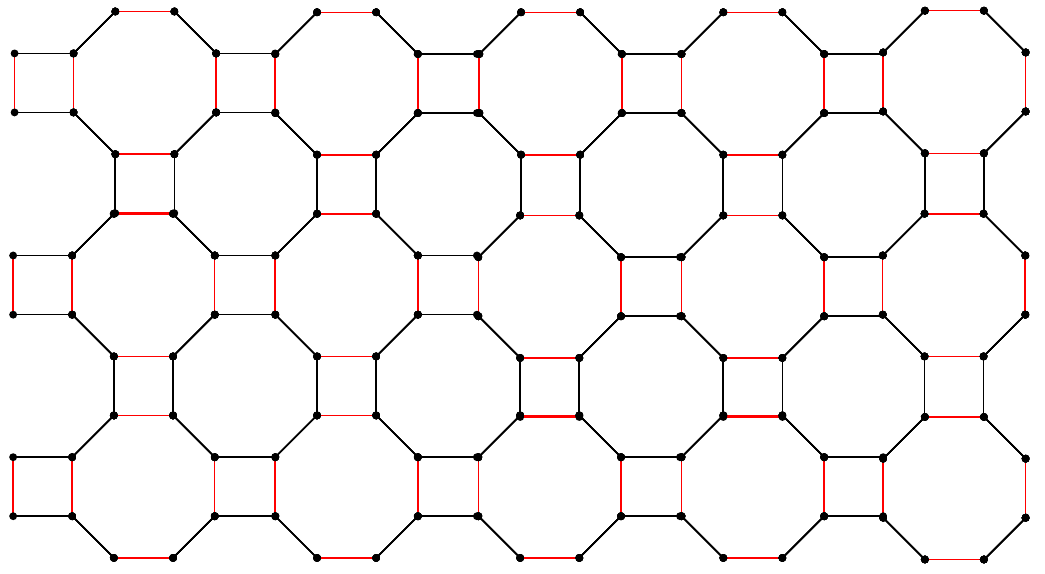}
    \caption{Perfect matching $M_0$ of $N_{5,3}$.}
    \label{sfig5}
  \end{minipage}\hfill
  \begin{minipage}{0.48\textwidth}
    \centering
    \includegraphics[width=\linewidth]{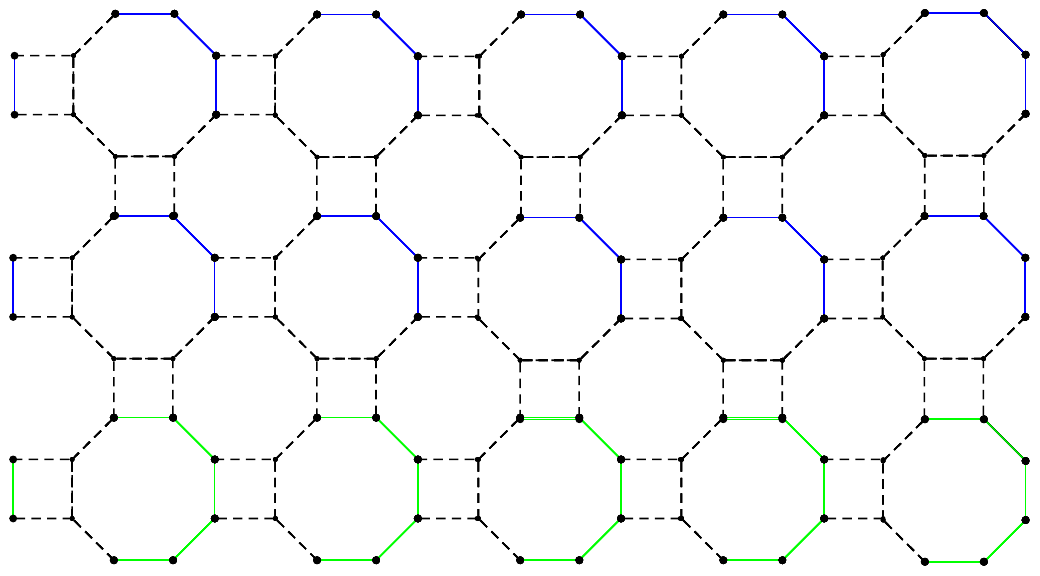}
    \caption{The subgraph \(N_{5,3}-V(S_0)\).}
    \label{sfig9}
  \end{minipage}
\end{figure}
Let
\[
S_0=\{a_{i,j}b_{i,j},\ c_{i,j}d_{i,j}\mid i\in\mathbb{Z}_p,\ j\in\mathbb{Z}_q\}
\setminus\{c_{i,0}d_{i,0}\mid i\in\mathbb{Z}_p\}
\]
be the subset of $M_0$ with $|S_0|=2pq-p=p(2q-1)$. For $i\in\mathbb{Z}_p$ and $j\in\mathbb{Z}_q$, the neighbors of $h_{i,j}$
lying in $V(S_0)$ are precisely $a_{i,j}$ and, whenever it exists,
$c_{i,j+1}$; those of $g_{i,j}$ are precisely $d_{i,j+1}$, whenever it
exists; and one neighbor $a_{i+1,j}$ of $f_{i,j}$ lies in $V(S_0)$, while
all remaining neighbors of $h_{i,j}$, $g_{i,j}$ and $f_{i,j}$ do not lie in $V(S_0)$. 
For $i\in\mathbb{Z}_p$ and $j\in\mathbb{Z}_{q,1}$, the neighbors
$b_{i+1,j}$ and $d_{i,j}$ of $e_{i,j}$ lie in $V(S_0)$, whereas $f_{i,j}$
does not. Hence $h_{i,j}g_{i,j}f_{i,j}e_{i,j}$ is an induced path of length
$3$ in $N_{p,q}-V(S_0)$. 
For $i\in\mathbb{Z}_p$ and $j=0$, the neighbors of $e_{i,0}$ and of
$c_{i,0}$ lying in $V(S_0)$ are precisely $b_{i+1,0}$ and $b_{i,0}$,
respectively, and all remaining neighbors of these two vertices do not lie in
$V(S_0)$. Hence $h_{i,0}g_{i,0}f_{i,0}e_{i,0}d_{i,0}c_{i,0}$ is an induced
path of length $5$ in $N_{p,q}-V(S_0)$. 
Therefore, \(N_{p,q}-V(S_0)\) is the disjoint union of the \(p(q-1)\) paths of length 3
\[
h_{i,j}g_{i,j}f_{i,j}e_{i,j}\qquad (i\in\mathbb{Z}_{p}, j\in\mathbb{Z}_{q,1}),
\]  indicated by the blue edges in Fig.~\ref{sfig9},
and the \(p\) paths of length 5
\[
h_{i,0}g_{i,0}f_{i,0}e_{i,0}d_{i,0}c_{i,0} \qquad (i\in\mathbb{Z}_{p}),
\] indicated by the green edges in Fig.~\ref{sfig9},
where the first family is empty when \(q=1\).
Since each such path has a unique perfect matching, so does \(N_{p,q}-V(S_0)\). Thus, \(S_0\) is a forcing set of \(M_0\), and therefore
\[
f(N_{p,q},M_0)\le |S_0|=p(2q-1).
\]
For the reverse inequality, every square of $N_{p,q}$ is an $M_0$-alternating cycle, and
$N_{p,q}$ contains $p(2q-1)$ pairwise disjoint squares. By Lemma~\ref{lem:alternating},
$f(N_{p,q},M_0)\ge p(2q-1)$. Together with the upper bound \(f(N_{p,q},M_0)\le p(2q-1)\), this gives $f(N_{p,q},M_0)=p(2q-1).$
Therefore,
\[
F(N_{p,q})\ge f(N_{p,q},M_0)=p(2q-1).
\]

It remains to prove that $F(N_{p,q})\le p(2q-1)$. Let $M$ be an  arbitrary perfect
matching of $N_{p,q}$. We choose an independent set $T=T'\cup T''$, where $$T'=\{e_{i,j},\ g_{i,j}\mid i\in\mathbb{Z}_p,\ j\in\mathbb{Z}_{q,1}\},$$
\[
T''=\{b_{i,0},\ d_{i,0},\ g_{i,0}\mid i\in\mathbb{Z}_{p,1}\}\cup\{a_{0,0},\ d_{0,0},\ g_{0,0}\},
\]
 indicated by the red vertices and green vertices, respectively, in Fig.~\ref{sfig10}. Then $$|T|=|T'|+|T''|=2p(q-1)+3(p-1)+3=p(2q+1).$$
 \begin{figure}[htbp]
  \centering
  \begin{minipage}{0.48\textwidth}
    \centering
    \includegraphics[width=\linewidth]{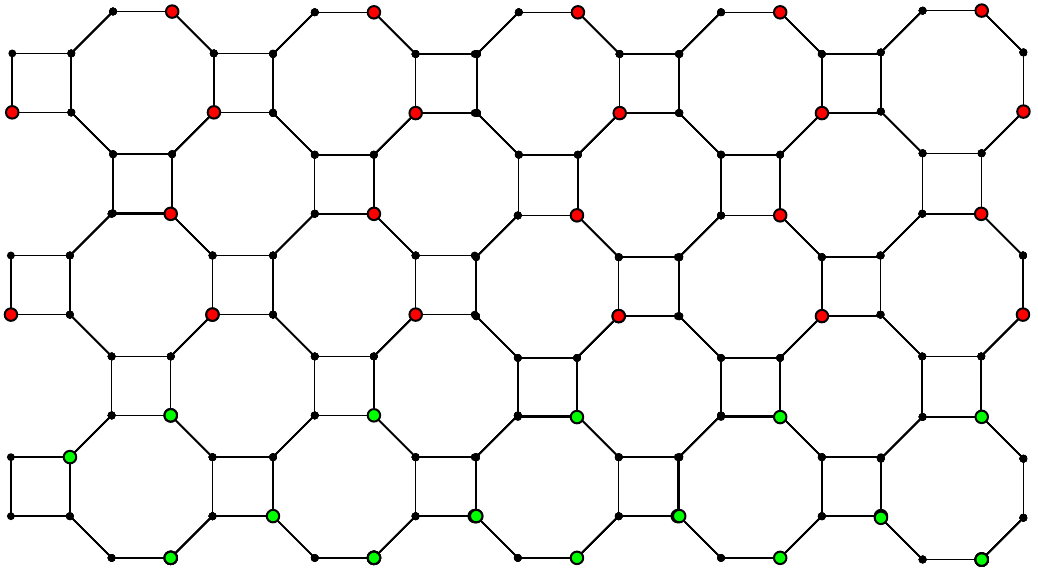}
    \caption{Marked vertices $T$ of \(N_{5,3}\).}
    \label{sfig10}
  \end{minipage}\hfill
  \begin{minipage}{0.48\textwidth}
    \centering
    \includegraphics[width=\linewidth]{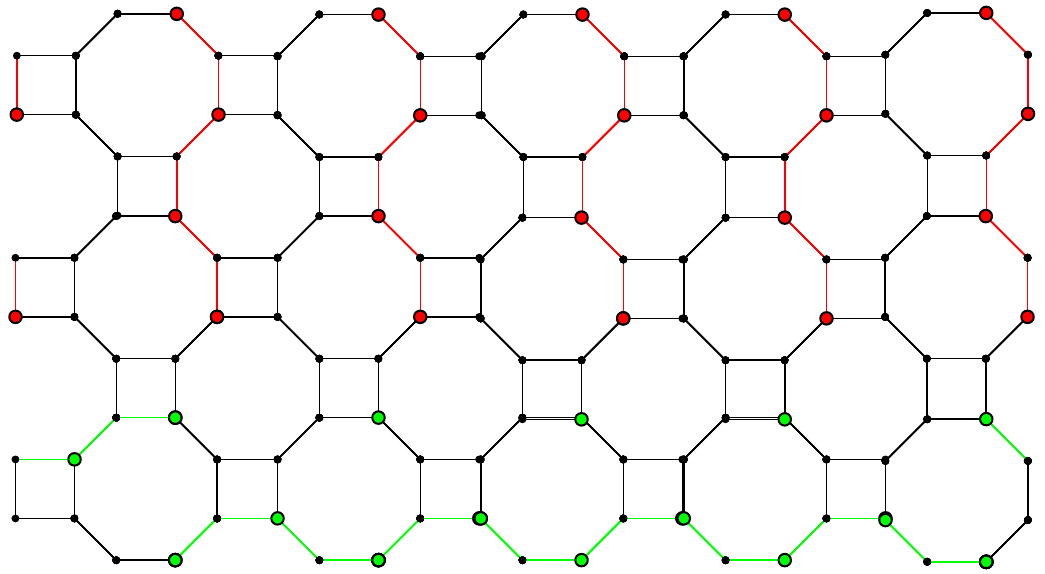}
    \caption{The union of 2-paths of \(N_{5,3}\).}
    \label{sfig11}
  \end{minipage}
\end{figure}
In particular, $T'=\emptyset$ for $q=1$, and $|T|=|T''|=3p$.
 Let $G'$ be the graph formed by the union of all 2-paths of $N_{p,q}$
whose endvertices lie in $T$. In what follows, we prove that $G'$  contains no
$M$-alternating cycle.

Note that a vertex belongs to $V(G')$ only if it lies in $T$ or has at least
two neighbors in $T$. It is clear that, the vertices \(a_{i,j}\), \(b_{i,j}\), and \(c_{i,j}\) are not in \(T\) for \(i\in\mathbb Z_p\) and \(j\in\mathbb Z_{q,1}\). We now show that they are not the  vertices of \(G'\).
For \(i\in\mathbb Z_p\) and \(j\in\mathbb Z_{q,1}\), since all three neighbors \(f_{i-1,j}\), \(b_{i,j}\) and \(h_{i,j}\) of \(a_{i,j}\) lie outside the independent set \(T\), we obtain \(a_{i,j}\notin V(G')\). For \(i\in\mathbb Z_p\) and \(j\in\mathbb Z_{q,1}\), we have \(e_{i-1,j}\in T\), while \(a_{i,j}\notin T, c_{i,j}\notin T.\) Thus \(b_{i,j}\) has only one neighbor in \(T\), and consequently \(b_{i,j}\notin V(G')\). For \(i\in\mathbb Z_p\) and \(j\in\mathbb Z_{q,1}\), all three neighbors \(b_{i,j}\), \(d_{i,j}\), and \(h_{i,j-1}\) of \(c_{i,j}\) lie outside \(T\), hence \(c_{i,j}\notin V(G')\).
For \(i\in\mathbb Z_{p}\) and \(j\in\mathbb Z_{q,1}\), \(h_{i,j}\notin T\). Its neighbor \(g_{i,j}\) lies in \(T\), while \(a_{i,j}\) and, for \(i\ne p\), \(c_{i,j+1}\) do not. Hence \(h_{i,j}\) has only one neighbor in \(T\), so \(h_{i,j}\notin V(G')\).
Therefore, the subgraph of \(G'\) consisting of the union of all 2-paths between the vertices of $T'$
is a disjoint union of \(p\) paths of length \(4(q-1)-2\)
\[e_{i,1}f_{i,1}g_{i,1}d_{i,2}e_{i,2}f_{i,2}g_{i,2}\cdots d_{i,q-1}e_{i,q-1}f_{i,q-1}g_{i,q-1} \qquad (i\in\mathbb{Z}_{p}).\]

The three vertices $b_{0,0}$, $c_{0,0}$ and $e_{p-1,0}$ are not in $T$, and
each of them has exactly one neighbor in $T$. Hence none of them is a vertex
of $G'$. For $i\in\mathbb{Z}_{p,1}$, since $h_{i,0}\notin T$ and $g_{i,0}$ is
its only neighbor in $T$, $h_{i,0}$ is not a vertex of $G'$. Similarly
$f_{i,0}$ for $i\in\mathbb{Z}_{p,-1}$ and $a_{i,0}$ for
$i\in\mathbb{Z}_{p,1}$ are not vertices of $G'$. Hence the subgraph of $G'$
formed by the union of all 2-paths between vertices of $T''$ is the disjoint
union of a path of length $4(p-1)$,
\[
  d_{0,0}e_{0,0}b_{1,0}c_{1,0}d_{1,0}e_{1,0}\cdots
  b_{p-1,0}c_{p-1,0}d_{p-1,0},
\]
of $p-2$ isolated vertices $g_{1,0},g_{2,0},\dots,g_{p-2,0}$, and of a path of
length four $g_{0,0}h_{0,0}a_{0,0}f_{p-1,0}g_{p-1,0}$.

 We now consider the union of all 2-paths with both endvertices in \(T\). When $q>1$, the \(p-2\) vertices $g_{i,0}$ are joined to the \(p-2\) paths of length \(4(q-1)-2\)
\[
e_{i,1}f_{i,1}g_{i,1}d_{i,2}e_{i,2}f_{i,2}g_{i,2}\cdots d_{i,q-1}e_{i,q-1}f_{i,q-1}g_{i,q-1}
\]
by a 2-path $g_{i,0}d_{i,1}e_{i,1}$, where \(i\in\mathbb{Z}_{p}\setminus \{0,p-1\}\), forming $p-2$ paths of length \(4(q-1)\); and  the remaining two paths of length \(4(q-1)-2\) are joined to the path $g_{0,0}h_{0,0}a_{0,0}f_{p-1,0}g_{p-1,0}$
of length \(4\) by the \(2\)-paths \(g_{0,0}d_{0,1}e_{0,1}\) and
\(g_{p-1,0}d_{p-1,1}e_{p-1,1}\), respectively, forming a path of length \(8q-4\).
Thus \(G'\) is a forest consisting of \(p\) disjoint paths, and therefore contains no \(M\)-alternating cycle.
 When $q=1$, the set $T'$ is empty
and $G'$ coincides with the union of the pieces from $T''$, which is already a forest. In either case $G'$ is a forest consisting of \(p\) pairwise disjoint
paths, so it contains no $M$-alternating cycle. Since $|M|=4pq$,
Lemma~\ref{lem:labeled} yields
\[
f(N_{p,q},M)\le |M|-|T|=4pq-p(2q+1)=p(2q-1).
\]
By the definition of the maximum forcing number, $F(N_{p,q})\le p(2q-1)$.
Together with $F(N_{p,q})\ge p(2q-1)$ proved above, we conclude that
$F(N_{p,q})=p(2q-1)$.
\end{proof}

\begin{thm}\label{thm:s-torus}
Let $p\ge 1$ and $q\ge 1$. Then the maximum forcing number
$F(H_{p,q})=2pq.$
\end{thm}

\begin{proof}
First, we prove that $F(H_{p,q})\ge 2pq$. Consider the perfect matching
\[
M_0=\{a_{i,j}b_{i,j},\ c_{i,j}d_{i,j},\ e_{i,j}f_{i,j},\ g_{i,j}h_{i,j}\mid i\in\mathbb{Z}_p,\ j\in\mathbb{Z}_q\}.
\]
of $H_{p,q}$, indicated by the red edges in Fig.~\ref{sfig12}.  Let
\[
S_0=\{a_{i,j}b_{i,j},\ c_{i,j}d_{i,j}\mid i\in\mathbb{Z}_p,\ j\in\mathbb{Z}_q\}
\]
be a subset of $M_0$. 
\begin{figure}[htbp]
  \centering
  \begin{minipage}{0.48\textwidth}
    \centering
    \includegraphics[width=\linewidth]{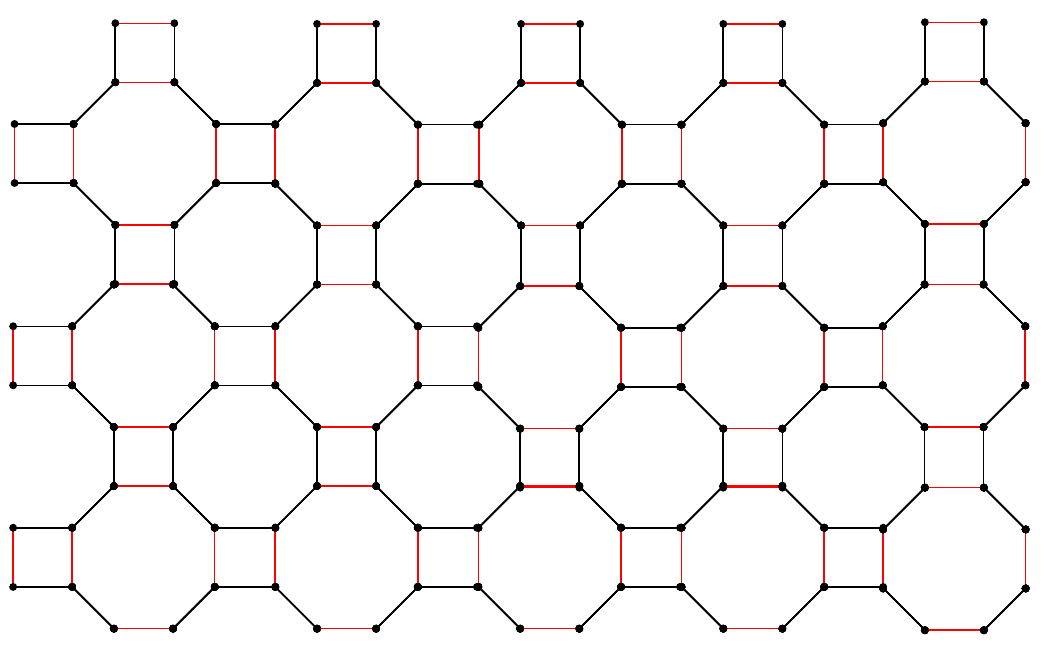}
    \caption{The perfect matching $M_0$ of \(H_{5,3}\).}
    \label{sfig12}
  \end{minipage}\hfill
  \begin{minipage}{0.48\textwidth}
    \centering
    \includegraphics[width=\linewidth]{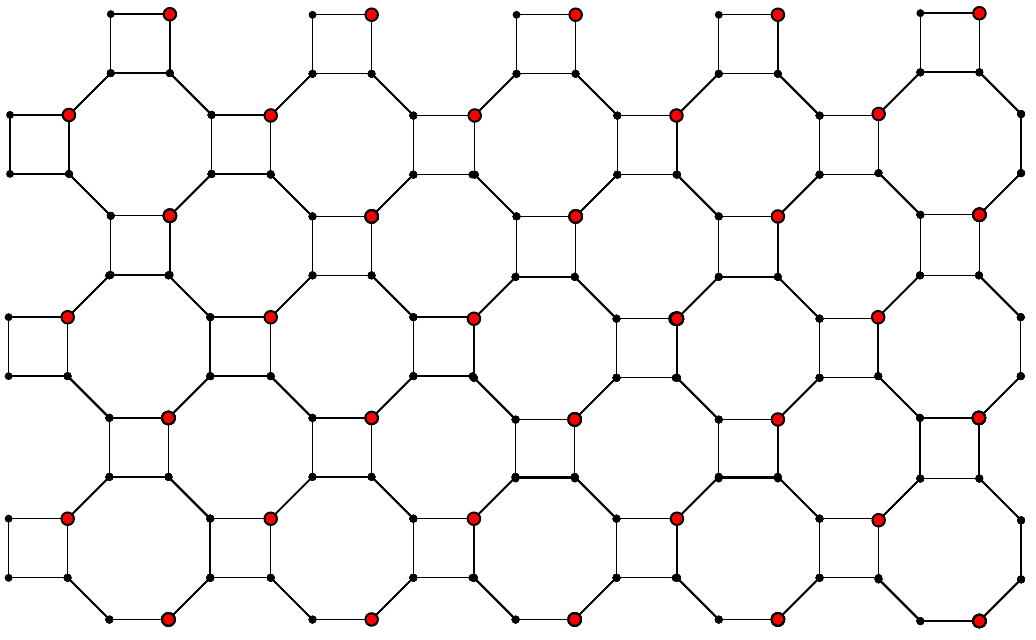}
    \caption{The marked vertices of \(H_{5,3}\).}
    \label{sfig11}
  \end{minipage}
\end{figure}

For $i\in\mathbb{Z}_p$ and $j\in\mathbb{Z}_q$, the neighbors of $h_{i,j}$
lying in $V(S_0)$ are precisely $a_{i,j}$ and $c_{i,j+1}$; those of $g_{i,j}$ are precisely $d_{i,j+1}$; one neighbor $a_{i+1,j}$ of $f_{i,j}$ lies in $V(S_0)$; the neighbors $b_{i+1,j}$ and $d_{i,j}$ of $e_{i,j}$ lie in $V(S_0)$; 
all remaining neighbors of $h_{i,j}$, $g_{i,j}$, $f_{i,j}$ and $e_{i,j}$ do not lie in $V(S_0)$. Hence $h_{i,j}g_{i,j}f_{i,j}e_{i,j}$ is an induced path of length
$3$ in $H_{p,q}-V(S_0)$. 
Since $H_{p,q}-V(S_0)$ consists of $pq$ pairwise disjoint paths of this form, and a path on an even number of vertices has a unique perfect matching, $H_{p,q}-V(S_0)$ has a unique perfect matching. Thus, $S_0$ is a forcing set of $M_0$, and $f(H_{p,q},M_0)\le |S_0|=2pq$. Since every square in $H_{p,q}$ is an $M_0$-alternating cycle and $H_{p,q}$ contains $2pq$ pairwise disjoint squares,  $f(H_{p,q},M_0)\ge 2pq$ by Lemma~\ref{lem:alternating}. Thus $f(H_{p,q},M_0)=2pq$, and so $F(H_{p,q})\ge f(H_{p,q},M_0)=2pq$ by the definition.

Next, we prove that $F(H_{p,q})\le 2pq$. Let $M$ be an arbitrary perfect matching of $H_{p,q}$ and take the marked vertices
\[
T=\{a_{i,j},\ d_{i,j}\mid i\in\mathbb{Z}_p,\ j\in\mathbb{Z}_q\}.
\]
Let $G'$ be the graph formed by the union of all 2-paths of $H_{p,q}$ having both endvertices in $T$. Since the distance between any two vertices in $T$ is at least \(3\), there is no 2-path joining two vertices of $T$. Hence $G'$ is the edgeless graph on the $2pq$ marked vertices. In particular $G'$ contains no $M$-alternating cycle. Noting that $|M|=4pq$ and $|T|=2pq$, Lemma \ref{lem:labeled} yields
\[
f(H_{p,q},M)\le |M|-|T|=4pq-2pq=2pq.
\]
By the definition of the maximum forcing number, $F(H_{p,q})\le 2pq$. Together with the inequality $F(H_{p,q})\ge 2pq$ proved above, we conclude that $F(H_{p,q})=2pq$, which completes the proof.
\end{proof}

\subsection*{Acknowledgment}
This work is supported by Natural Science Foundation of Inner Mongolia Autonomous Region (Grant No. 2025QN01008), and first-class discipline research special project (Grant No. YLXKZX-NGD-055).

%%\noindent{\normalsize \textbf{Author Contributions}} The authors have the same contributions to the paper. All authors reviewed the manuscript.
%%\vspace{3mm}

%\vspace{3mm}
%\noindent{\normalsize \textbf{Data availability}}
%No data was used for the research described in the article.

%\vspace{3mm}
%\noindent{\normalsize \textbf{Declarations}}
%\vspace{3mm}

%\texttt{Conflict of interest} The authors declare no conflict of interests.


\begin{thebibliography}{99}
\setlength{\itemsep}{-.2mm}


\bibitem{afshani04}
P. Afshani, H. Hatami, E.~S. Mahmoodian, On the spectrum of the forced matching number of graphs,
\emph{Australas. J. Combin.} \textbf{30} (2004), 147--160.

\bibitem{adams04}
P. Adams, M. Mahdian, E.~S. Mahmoodian, On the forced matching numbers of bipartite graphs,
\emph{Discrete Math.} \textbf{281} (2004), 1--12.

\bibitem{diudea05}
M.~V. Diudea, Covering nanostructures,
in: \emph{Nanostructures---Novel Architecture}, NOVA, New York, 2005, pp.~203--242.



\bibitem{harary91}
F.~Harary, D.~J. Klein, T.~P. \v{Z}ivkovi\'{c}, Graphical properties of polyhexes: perfect matching vector and forcing,
\emph{J. Math. Chem.} \textbf{6}(1) (1991), 295--306.

\bibitem{iranmanesh07}
A. Iranmanesh, B. Soleimani, PI index of $TUC_4C_8(R)$ nanotubes,
\emph{MATCH Commun. Math. Comput. Chem.} \textbf{57} (2007), 251--262.

\bibitem{jiangzhang16}
X. Jiang, H. Zhang, The maximum forcing number of cylindrical grid, toroidal 4--8 lattice and Klein bottle 4--8 lattice,
\emph{J. Math. Chem.} \textbf{54}(1) (2016), 18--32.

\bibitem{klein87}
D.~J. Klein, M.~Randi\'{c}, Innate degree of freedom of a graph,
\emph{J. Comput. Chem.} \textbf{8}(4) (1987), 516--521.

\bibitem{kleinerman06}
S. Kleinerman, Bounds on the forcing numbers of bipartite graphs,
\emph{Discrete Math.} \textbf{306}(1) (2006), 66--73.


\bibitem{stefu04}
M. Stefu, M.~V. Diudea, Wiener index of $C_4C_8$ nanotubes,
\emph{MATCH Commun. Math. Comput. Chem.} \textbf{50} (2004), 133--144.

\bibitem{xudeng07}
L. Xu, H. Deng, PI indices of tori $T_{p,q}[C_4,C_8]$ covered by $C_4$ and $C_8$,
\emph{MATCH Commun. Math. Comput. Chem.} \textbf{57} (2007), 485--502.

\bibitem{xuzhang09}
X. Xu, H. Zhang, Hosoya polynomials of $TUC_4C_8(S)$ nanotubes,
\emph{J. Math. Chem.} \textbf{45} (2009), 488--502.


\bibitem{zhang25}
Y. Zhang, X. He, Q. Liu, Forcing, anti-forcing, global forcing and complete forcing on perfect matchings of graphs---A survey,
\emph{Discrete Appl. Math.} \textbf{376} (2025), 318--347.
\end{thebibliography}
\end{document}